\input amstex
\documentstyle{amsppt}

\def\L{{\Cal L}}
\def\Lat{\Lambda }
\def\Lv{\Lambda_v }
\def\Q{\Bbb Q}
\def\R{\Bbb R}
\def\C{\Bbb C}

\def\Kstar{{K}^{\vphantom{l}{\times}}}

\def\Kvn{(K_v)^n}
\def\ZK{{\Cal O}_{K}}
\def\Ov{{\Cal O}_{v}}


\def\dx#1{d_{\X}(#1)}

\def\hop{{H}^{\operatorname{op}}}
\def\hopf{{H}^{\operatorname{op}}_{\F}}

\def\sph{{H}_{s}}
\def\M#1#2{{\Cal M}^{\sssize{#1}}_{\sssize{#2}}}

\def\abv#1{|#1|_{{}_{\ssize v}}}

\def\F{\Cal F}
\def\E{\Cal E}


\def\Mn{\operatorname{M}_n}

\def\Endk{\operatorname{End}(K^n)}
\def\End#1{\operatorname{End}\bigl(#1\bigr)}

\def\AGL#1{\operatorname{GL}_{#1}(K_{\Bbb A})}
\def\GLnk{\operatorname{GL}_n(K)}

\def\projkn{{\Bbb P} \bigl(\Endk\bigr)}
\def\om{\Omega\bigl({\Bbb P} \bigl(\Endk\bigr), B\bigr)}
\def\proj#1{{\Bbb P} \bigl(\operatorname{End}(#1)\bigr)}
\def\P{\projkn}
\def\Pn{{\Bbb P}^n}


\def\normm#1{\|\,#1\,\|}
\def\normdv#1{\|\,#1\,\|_{{}_{\ssize v}}^{d_v}}
\def\nv#1{\|\,#1\,\|_{{}_{\ssize v}}}
\def\nvb#1{\|\,#1\,\|_{{}_{\ssize b,v}}}


\def\tensork{{\otimes}_{{}_{K}}}
\def\Sp#1{\operatorname{sp}(#1)}
\def\rov#1{{\rho}_{{}_{\ssize v}}(#1)}
\def\rovdv#1{{\rho}_{{}_{\ssize v}}(#1)^{d_v}}


\def\ker#1{\operatorname{ker}#1}

\def\rank#1{\operatorname{rank}\left(#1\right)}

\def\dim#1#2{\operatorname{dim}_{\sssize #1}#2}

\def\mz{-\{\zero\}}
\def\fine{\quad{\ssize \blacksquare}}
\def\res#1{\lower2pt\hbox{${|}_{#1}$}}
\def\iso{\simeq}

\def\st{\,\bigm|\,}

\def\cs{\Cal P}

\def\ra{\longrightarrow}
\def\lqtwo{{\ell}^{@!@!{\raise1.5pt\hbox{$\sssize 2$}}}}


\def\X{\Cal X}
\def\Y{\Cal Y}

\def\E{\Cal E}

\def\nxv#1{\|\,#1\,\|_{{}_{\ssize \X_v}}}
\def\nxdv#1{\|\,#1\,\|_{{}_{\ssize \X_v}}^{d_v}}
\def\x{\bold x}
\def\y{\bold y}
\def\z{\bold z}
\def\e{\bold e}
\def\t{\bold t}

\def\pl{{\Cal P}_\ell}
\def\zero{\pmb 0}
\def\b{\underline{b}}

\def\cl#1{[#1]}

\def\ba{\setminus}
\def\comp{{\ssize \circ}}

\topmatter
\title A Gelfand-Beurling type formula for heights on endomorphism rings
\endtitle
\author Valerio Talamanca
\endauthor
\affil{Universiteit van Amsterdam}\endaffil
\address{KdV Instituut voor Wiskunde, Plantage Muidergracht 24, 1018 TV Amsterdam, 
The Netherlands}
\endaddress
\abstract{Let $V$ be a finite dimensional vector space over a number field.
In this paper we prove a limit formula for heights on the endomorphism ring of $V$, 
which can be considered as the analogue for both the
Gelfand-Beurling formula for the
spectral radius on a Banach Algebra and Tate's averaging
procedure for constructing canonical heights on Abelian Varieties. We also
prove a version of Northcott's finiteness theorem.}
\endabstract
\email{valerio\@wins.uva.nl}
\endemail
\date{October 4, 1999}
\enddate
\rightheadtext{A Gelfand-Beurling type formula for heights}
\thanks{Research partially supported by a post-doc fellowship
of the University of Padova and C.N.R. The author would also like to thank
the Mathematics Department of Universit\`a di Roma III for its kind hospitality.}
\endthanks
\endtopmatter
\document
\head Introduction\endhead 
 Let $A$ be an abelian variety defined over a number field $K$. 
Fix an ample and symmetric invertible sheaf $\Cal L$ on $A$. Suppose that 
$\phi:A\rightarrow \Pn$ is an injective morphism associated to $\L$, i.e.
$\phi^*{\Cal O}_{\Pn}(1)\iso\L$. Let $H$ be the Northcott-Weil $\ell^2$-height on 
$\Pn ({\overline K})$ (we will recall the definition of $H$ in section 1), and 
set $h_{\phi}=\log (H\comp\phi):A({\overline K})\rightarrow \R$. The function $h_{\phi}$ is called 
a (logarithmic) Northcott-Weil height on $A$ associated to $\L$. Clearly $h_{\phi}$ is not  uniquely determined by $\L$, but it can be shown (see e.g. \cite{9}) that the various functions constructed using  different choices of $\phi$ (always subject to the condition $\phi^*{\Cal O}_{\Pn}(1)\iso\L$) all lie in the same class modulo bounded functions. A lemma of J. Tate, as stated in \cite{8, pp. 29-30},
allows us to choose, amongst all the Northcott-Weil heights
associated to $\Cal L$, a canonical one having good functorial properties. If we denote by  $[n]:A
\ra A$ the multiplication by 
$n$ map, then we can  
explicitly compute the canonical height $\hat{h}_{\Cal L}$ via Tate's averaging 
procedure:
$$
\hat{h}_{\Cal L}(P)=\lim_{k\to\infty}\frac{h_{\phi}([n^k]P)}
{n^{2k}}.\tag{0.1}
$$
Let us stress that the function $\hat{h}_{\Cal L}$ is independent of the choice of $\phi$ and $[n]$.\par
The purpose of this paper is to prove an analogous formula for 
heights on the endomorphism rings of a finite dimensional $K$-vector space. 
As for the case of abelian varieties one needs some additional structure in order 
to be able to define heights. In this paper we work with heights associated to 
an adelic norm on $V$. We will provide all the relevant definitions in section 1. 
Let $V$ be a finite dimensional $K$-vector space and $\End{V}$ its endomorphism ring.
To any adelic norm $\F$ on $V$ we associate an height function $\hopf $ on $\End{V}$. 
On $\End{V}$ there also exists another height function, $H_s$,  which is called the spectral 
height, and is substantially different from the operator heights being 
defined as the product of the local spectral radii.
In our situation the  spectral height will play the role of the canonical height. 
In fact
our first main result is:
\proclaim{Theorem A}
Let $V$ be a finite dimensional $K$-vector space, and $\F$ an adelic norm  
on $V$. 
Then, for all $T{\in}\End{V}$, we have
$$
\lim_{k\to\infty} \hopf (T^k)^{{}^{\ssize \frac1k}}=\sph (T).\tag{0.2}
$$
\endproclaim
Clearly (0.2)
 is the analogue  of Tate's averaging procedure (0.1) in this setting.
Recall that the  Gelfand-Beurling 
formula for the spectral radius on a 
complex Banach algebra with $1$, $\bigl(A, \|\cdot\|\bigr)$, states that  
$$\dsize\lim_{k\to\infty}{\normm{a^k}}^{\frac1k}=\rho
(a), $$
where $\rho (a)=\sup_{\lambda{\in}\Sp a}|\lambda |$ and 
 $\Sp a=\{\lambda{\in}\C\st a-\lambda\cdot 1 \text{ is not invertible}\}$.
Therefore (0.2) can also be considered as the 
global analogue of the (finite dimensional) Gelfand-Beurling 
formula. 
Let us remark that the Gelfand-Beurling formula, 
and its real and $p$-adic counterparts, are actually used in the proof of theorem A. \par  
 Let $\projkn$ denote the projective space
associated to $\Endk$. The operator height, being homogeneous,
descends to a real valued function on $\projkn$. The other significant result that we will present is the extension to this setting of
 Northcott's theorem stating the finiteness of the set of points of bounded heights in 
projective spaces:
\proclaim{Theorem B} Let $V$ be a finite dimensional $K$-vector space and 
$\F$ a regular  adelic norm on $V$. Then the set
$$
\bigl\{[T]{\in}\proj{V}\st \rank{T}\geq 2 \text{\,and\,}
\hopf (T)\leq B\bigr\}
$$
is finite for every $B\geq 1$. 
\endproclaim
It is necessary to exclude rank one transformations from the above statement,
see section 4 for more details.
\par 

The paper is organized as follows: in section 1 we set our notations
and we introduce the height functions that we will use. In section 2 we prove a reduction lemma
which shows that is sufficient to prove our main result for the $\ell^2$
height on $K^n$. Section 3 contains our main technical result (theorem 3.3) which is a comparison 
result between different heights on $\End{K^n}$.
Section 4 is devoted to the proof of the main results. 
\vskip.1in {\it Acknowledgements} I would like to thank D. Masser and  P. Stevenhagen for many comments and suggestions. It is a pleasure to express my gratitude to Alan Mayer for encouragement and advice.
\vskip.2in
\head  1. Heights\endhead
Let $K$ be a number field of degree $d$ over $\Q$. We denote by
$\M{}{K}$ the set of equivalence classes of absolute values of $K$; by
$\M{0}{K}$ (respectively $\M{\infty}{K}$) the subset of $\M{}{K}$ consisting of 
the equivalence classes of non-archimedean (resp. archimedean) absolute values. 
If $v{\in}\M{0}{K}$, $v|p$, we normalize $\abv{\cdot}$ 
by requiring that $\abv{p}=p^{-1}$; while if  $v{\in}\M{\infty}{K}$ we then 
normalize $\abv{\cdot}$ by requiring that its restriction to $\Q$ is the 
standard archimedean absolute value.
Let $K_v$ be the completion of $K$ with respect to ${\abv{\cdot}}$. We denote by $n_v$ the local 
degree, and set $d_v=n_v/d$. \par
Before giving the definition of adelic norm we need to recall a few facts about
lattices over number field and their completions, see \cite{13, ch.2 \& 5} for more details.
Given $v{\in}\M{0}{K}$ we denote by $\Ov$ the closure of $\ZK$ (the ring of integers 
of $K$) in $K_v$. A $K_v$-lattice in a finite dimensional $K_v$-vector space 
is a compact and open $\Ov$-module. Let $W$ be a $K_v$-vector space and 
let $M\subset W$ be a $K_v$-lattice $M$. The norm 
associated to $M$, $N_M:W\rightarrow \R$, 
is defined by
$$
N_M({\x})=\inf_{\gamma{\in}K_v^\times, \gamma\x{\in}M}\abv{\gamma}^{-1}.
$$
Let $V$ be an $n$-dimensional $K$-vector space.
An $\ZK$-module $\Lat$ in $V$ is called a $K$-lattice if it is finitely generated and contains a basis 
of $V$ over $K$. Given $v{\in}\M{0}{K}$ we denote by $\Lv$  the closure of 
$\Lat$ in $K_v^n$. 

\par
By an {\it adelic norm}\footnote{Some 
comments on why we have chosen this definition of adelic norm will be provided
in the remark following the definition of the operator height associated to $\F$.} (cf. \cite{14}) on $V$ we mean a collection
$\F=\{N_v\,,\, v{\in}\M{}{K}\}$ of norms $N_v:V_v=V\tensork K_v\longrightarrow \R
$, having the following properties:
\roster
\item"(a)" $N_v$ is a norm with respect to $\abv{\cdot}$. Moreover, 
if $v{\in}\M{0}{K}$, then $N_v$ is ultrametric, i.e.
$N_v(\x+\y){\leq} \max\{N_v(\x),N_v(\y)\}$.
\item"(b)" There exists a $K$-lattice $\Lat$, such that $N_v$ is the norm associated
to $\Lv$ for all but finitely many $v{\in}\M{0}{K}$.
\endroster
A moment of reflection shows that if  $\Lambda$ is a $K$-lattice in $V$, then there exists a basis 
$\{\y_1,\ldots ,\y_n\} $ such that for all but finitely many $v{\in}\M{0}{K}$, $\Lv$ is the 
$\Ov$-module generated by $\{\y_1,\ldots ,\y_n\}$. In particular for each $\x{\in}V$
the set 
$\{v{\in}\M{0}{K}| N_{\Lv}(\x)\neq 1\}
$
is finite. Therefore, given an adelic norm $\F$ on $V$, it makes sense to set 
$$
 H_{\F}(\x)=\underset{v{\in}\M{}{K}}\to\prod N_v(\x)^{d_v},
$$
for all $\zero\neq\x{\in}V$. We set by definition $H_{\F}(\zero)=1$. The 
function $H_{\F}:V\rightarrow \R$ so defined is called the {\it height}  associated to $\F$.
Note that the product formula implies that $H_{\F}$ is homogeneous, i.e.
$H_{\F}(\lambda \x)=H_{\F}(\x)$ for all $\lambda{\in}\Kstar$. 
\par

\example{Examples}\,(a)\quad  Let $V=K^n$.
Set
$$
 \nv {\x}=\cases \left
(\sum\limits_{i=1}^n |{x_i}|_{{}_{\ssize v}}^2 \right)^{1\over 2}  & \quad
\text {if $ v{\in}\M{\infty}{K} $}\\ \sup_{1\leq i\leq n} \abv{x_i} &\quad
\text {if $v{\in}\M{0}{K}$. } \endcases  
$$  
Then ${\Cal E}=\{\nv{\cdot}\}_{v{\in}\M{}{K}}$ is an adelic norm
on $K^n$ and its associated height function $H_{\Cal E}$ is the
$\ell^2$ Northcott-Weil height. 
By changing the $\ell^2$-norms at the archimedean places 
into either $\ell^1$ or $\ell^{\infty}$-norms
one recovers the other two Northcott-Weil heights that are
commonly used. Note that $H_{\E}$ is invariant under field extensions. 
\par
\noindent (b)\quad Let $V$ be an $n$-dimensional $K$-vector space and
$\b=\{\y_1,\ldots ,\y_n\}$ a basis of $V$ over
$K$. Let $\iota_{\b}$ be the isomorphism of $V$ to $K^n$
defined by mapping $\y_i$ to $\e_i$, where $\{\e_1,\ldots ,\e_n\}$ is the
canonical basis of $K^n$. Set $\nvb{\x}=\nv{\iota_{\b}(\x)}$, then
${\Cal F}_{\b}=\{\nvb{\cdot}, v{\in}\M{}{K}\}$ is an
adelic norm.
\par\noindent (c)\quad Let ${\Cal T}=(T_v)$ be an element of $\AGL{n}$, the adele group
of $\GLnk$. Define $N_{v} :K_v^n \rightarrow \R$ by setting
$
N_v (\x)=\nv{T_v(\x)}.
$
Then ${\Cal F}_{\Cal T} =\{N_v,\,v{\in}\M{}{K}\}$ is an  adelic norm on $K^n$. The  associated height
$H_{{\Cal F}_ {\Cal T}}$ was used  by D.Roy and J.Thunder in \cite{6}, where it is called the twisted height associated to ${\Cal T}$.
\endexample
Let us point out that
the height arising from the adelic norms defined in the examples include all the
height functions that are commonly used in the literature. \par
Let $\F$ be an adelic norm on $V$. The function $\hop_{\F}:\End{V}\rightarrow \R$ 
defined by setting
$$
\hopf (T)=\sup_{\x{\in}V}\frac{H_{\F}\bigl(T(\x)\bigr)}{H_{\F} (\x)}
$$
is called the {\it operator height} associated to $\F$ (or to $H_{\F}$). 
The following properties of $\hopf$ are an immediate consequence of the above
definition:
\roster
\item $\hopf (\lambda T)=\hopf (T)$
\item $\hopf (TS)\leq \hopf (T)\hopf (S)$.
\endroster
\par
Note that property (1) ensures us that $\hopf$ descends to a well-defined function
on  $\P$. We will see in the next section that $\hopf$ is well defined 
meaning that the $\hopf (T)<\infty$. 
\remark{Remark}  
Condition (b) in the definition of adelic norm is in some sense a strong one. In fact it is
immediate to verify that it is equivalent to require that
there exists a basis $\b$ of $V$ such that  $N_v =\nvb{\cdot}$ for all but finitely many
$v{\notin}\M{0}{K}$. A possible and natural way of relaxing condition (b) would be to
require only that 
the set $\{v{\in}\M{}{K}\st N_v(\x)\neq 1\}$ is finite for all $0\neq\x{\in}V$.
This condition is indeed sufficient to have a well-defined height function on $V$ attached to 
$\F$, but not for ensuring that the operator height is well defined, as shown
by the following example. Consider the family of norms 
$\F=\{N_p \st \M{}{\Q}\}$, where
$
N_p:\Q_p^2 \rightarrow \R
$
is defined as
$$
N_p(x_1, x_2)=\cases
\max\{|{x_1}|_p, |px_2|_p\}  & \quad
\text {if $ p\neq \infty $}\\ \max\{|{x_1}|, |x_2|\} &\quad
\text {if $p=\infty$. } \endcases  
$$  
Clearly $N_p$ is an ultrametric norm for any $p\neq \infty$. Moreover
given $(x_1, x_2)$ there are only finitely many $p$'s for which $N_p(x_1,x_2)\neq 1$.
Let $q$ be a prime
$$
H_{\F}(q, 1)=q\cdot N_q(q,1)=qq^{-1}=1,\tag{1.1}
$$
while
$$
H_{\E}(q,1)=\sqrt{q^2+1}.
$$
It follows that
$$
\frac{H_{\E}(q,1)}{H_{\F}(q, 1)}=\sqrt{q^2 +1}
$$
which is clearly unbounded as $q\to \infty$. This fact contradicts the conclusion of lemma 2.1 
(see next section). Therefore $\F$ cannot be 
an adelic norm. Let $T=\left(\smallmatrix 1&-1\\ 0&1\endsmallmatrix\right)$. We want to show that 
$$
\sup_{\x{\in}\Q^2}\frac{H_{\F}\bigl(T(x)\bigr)}{H_{\F}(\x)}=+\infty.
$$
We will accomplish this by exhibiting a sequence $\{\x_n\}$ of vectors in $\Q^2$ 
such that $H_{\F}(\x_n)=1$ and $H_{\F}\bigl(T(\x_n)\bigr)\to\infty$ as $n\to\infty$. 
First of all note that if $x=(2^mq,1)$, with $q$ being any integer, then
$
H_{\F}(x)\geq 2^{m-1}$. Then, let  $\{p_n\}$ be a sequence of prime numbers such that 
 $p_n\equiv 1$ mod $2^{n}$. The existence of such a sequence is guaranteed by repeated 
applications of Dirichlet's theorem on primes in arithmetic progressions. Then let 
$\{\x_n=(p_n, 1)\}$. By (1.1) $H_{\F}(\x_n)=1$, on the other hand $T(\x_n)=(p_n-1, 1)=(2^nq_n,1)$ and so
 $H_{\F}\bigl(T(x_n)\bigr)\geq 2^{n-1}$.
\endremark
\vskip.1in
We conclude this section by giving the definition of the spectral height. 
Let us recall
the definition of the local spectral radii.
Let $F$ be a complete local field and $W$ a finite dimensional
$F$-vector space. The {\it spectral radius} of $T{\in}\End{W}$ is
$
\rho_F(T)=\sup_{\lambda{\in}\Sp{T}}|\lambda|_{F(\lambda)},
$
where $\Sp{T}$ is the set of characteristic roots of $T$, and
$|\cdot|_{F(\lambda)}$ is the unique extension of
$|\cdot|_F$ to $F(\lambda)$.
 Given $T{\in}\End{V}$ we set
$
\rho_v(T)=\rho_{K_v}(T_v),
$
where $T_v$ is the extension of $T$ to $V_v$ by $K_v$-linearity. If $T$ is not nilpotent we set
$$
\sph (T)=\underset{v{\in}\M{}{K}}\to\prod \rovdv{T}.
$$
We set $\sph (T)=1$ for any nilpotent transformation. 
The function thereby  defined is called 
 the {\it spectral height} and  enjoys the following
properties:
{\it \roster
\item"(S1)" $\sph (\lambda T)=\sph (T) .$
\item"(S2)"  $\sph (T)\geq 1 .$
\item"(S3)" $\sph (T^k)=\sph (T)^k$.
\item"(S4)" $\sph$ is invariant under conjugation.
\item"(S5)" If $T_s$ is the semisimple, part of $T$ then
$
\sph (T)=\sph(T_s).
$
\item"(S6)" If $T,T^{\sssize \prime}{\in}\End{V}$ commute, $
\sph (TT^{\sssize \prime})\leq\sph (T)\sph (T^{\sssize \prime})$.
\item"(S7)" $\sph$ is invariant under field extension.
\endroster}
Properties (S3)-(S6) are direct consequences of the behavior of the
spectrum  under the various
operations considered (see \cite{5}). Property (S1) follows from the
product formula while (S7) is derived in a standard way
from the formula for local degrees
(see \cite{5, ch.3 \S 1}). Finally
 (S2) follows from (S1) and (S7).
\vskip.2in
\head  2. A reduction lemma\endhead
The main goal of this section is to show, that in order to establish
theorems A and B
in full generality, it is sufficient to prove them
for $\hop$ on $K^n$ (recall that $\hop$ is the operator height
associated to the standard $\ell^2$ adelic norm on $K$).
\proclaim{Lemma 2.1} {\it Let $\F=\{N_v, v{\in}\M{}{K}\}$ be an adelic norm on
$V$. Then there exists a constant $C>0$ and an isomorphism
$\iota : V\rightarrow K^n$ such that
$$
C^{-1}\hop (T^{\iota})\leq  H^{op}_{\Cal F} (T)
\leq C\hop (T^{\iota}),
$$
where $T^{\iota}=\iota\circ T\circ \iota^{-1}{\in}\Endk$.}
\endproclaim
\demo{Proof} Let $\Lambda$ be a $K$-lattice such that $N_v=N_{\Lv}$ for all but finitely many
$v{\in}\M{0}{K}$. By definition $\Lambda$ contains a basis  $\b=\{\y_1,\ldots , \y_n\}$
of $V$ over $K$. It follows that there exists a finite set of absolute values $\cs$, containing $\M{\infty},
{K}$, such that $N_v =\nvb{\cdot}$ for all 
$v{\notin}\cs$.
Let $\iota : V\rightarrow K^n$,  $\y_i\mapsto\e_i$. If $v{\notin}{\Cal S}$, then
$N_v(\x)=\nvb{\x}=\nv{\iota(\x)}$. On the other hand 
 all the norms on a finite dimensional vector space over a
complete field are equivalent, thus $\cs$ being finite we find that there exists 
$C_1, C_2>0,$ such that 
$$
C_1\cdot
H\bigl(\iota (\x)\bigr)\leq H_{\F} (\x)\leq C_2\cdot H\bigl(\iota (\x)\bigr).
$$
Then
an easy computation shows 
$$
C^{-1}\hop (T^{\iota})\leq  \hopf (T)
\leq C\hop (T^{\iota})
$$
with $C=C_2/C_1 >0.\fine$
\enddemo
 Let $\F$ be a  adelic norm on $V$. Given $B>0$, set
$$\Omega\bigl(\proj{V}, \hopf, B\bigr)=\bigl\{[T]{\in}\proj{V}\st \rank{T}\geq 2 \text{\,and\,}
\hopf (T)\leq B\bigr\}.
$$ 
The next lemma achieves the goal of this section.
\proclaim{Lemma 2.2} {\it If theorems A and B
hold for $\hop$ on $K^n$, then they hold
for the operator height associated to any  adelic norm
on an $n$-dimensional $K$-vector space.}
\endproclaim
\demo{Proof} 
Let  $\iota : V\rightarrow K^n$ and $C>0$ be as in the conclusion
of the previous lemma.
It follows that for every $B>0$, the map $T\mapsto T^{\iota}$ gives rise to an injection
$
\Omega\bigl(\proj{V}, \hopf, B\bigr)\hookrightarrow
\Omega\bigl(\projkn, \hop, BC\bigr).
$\par 
Regarding the assertion of theorem A we first note that 
if theorem A holds for $\hop$ on $K^n$, then 
$$
\limsup_{k\to\infty}\hopf
(T^k)^{\frac1k}{\leq}\limsup_{k\to\infty} \bigl(C\hop(({T^{\iota}})^k)\bigr)^{\frac1k}
=\lim_{k\to\infty} 
\bigl(\hop((T^{\iota})^k)\bigr)^{\frac1k}=\sph (T^{\iota})$$
and 
$$
\liminf_{k\to\infty}
\bigl(C^{-1}\hop({T^{\iota}}^k)\bigr)^{\frac1k}{\leq}\liminf_{k\to\infty}
\hopf (T^k)^{\frac1k}=\lim_{k\to\infty} 
\bigl(\hop((T^{\iota})^k)\bigr)^{\frac1k}=\sph (T^{\iota}).
$$
Combining these inequalities with the fact that $\sph (T)=\sph (T^{\iota})$ yields the desired limit
formula for $\hopf. \fine$
\enddemo
Thanks to lemma 2.2. we can restrict our attention to $H$ on $K^n$. 
Let us introduce the heights, related to $H$, that will be useful in the sequel. 
If
$\X$ is one-dimensional we set $H(\X)=H(\x)$, where $\zero\neq\x{\in}\X$ is any
non-zero element. A subspace
$\X\subset K^n$ of dimension $\ell$ determines a one-dimensional subspace $\pl
(\X)$ of $K^m$, where $m=\binom n\ell$ (via the Pl\"ucker map), allowing us to set
$H(\X)=H\bigl(\pl (\X)\bigr)$, cf. \cite{7} where this height was firstly  introduced.
\par Next we want to introduce an auxiliary height function on $\End{K^n}$, 
which is defined as the product of local operator norms. More precisely
the norm $\|\cdot\|_v$
induces a norm on $\End{\Kvn}$, defined by $$ \nv{S}=\sup_{\x{\in}\Kvn\mz}
\frac{\nv{S(\x)}}{\nv{\x}}. $$
Explicitly, given $S=(s_{ij}){\in}M_n(K_v)$, we have  $\nv{S}=\sup_{1\leq i,j \leq n} \abv{s_{ij}}$ if $ v{\in}\M{0}{K}$. \linebreak
If  $v{\in}\M{\infty}{K}$, then  $
\nv{T}=\sup_{\lambda{\in}\Sp{T^*T}} \sqrt{\lambda}
$
where $T^*$ is the adjoint of $T$.
It follows that for $0\neq T{\in}\End{K^n}$, we have $\nv{T}=1$ for all but finitely many $v$'s. Therefore 
it makes sense to set:
$$
H(T)=\underset
{v{\in}\M {}{K}}\to{\prod}
\normdv{T}.
$$
As usual, we set $H(0)=1$. The elementary
properties of $H:\End{K^n}\rightarrow \R$ are:
\roster
\item"(H1)\," $H(\lambda T)=H(T)$;
\item"(H2)\," $H(T)\geq 1$; 
\item"(H3)\," $H(TT')\leq H(T)\cdot H(T')$; 
\item"(H4)\," $\hop (T)\leq H(T)$
\item"(H5)\," {\it Let $\phi:\Endk\rightarrow K^{n^2}$ be an isomorphism which 
assigns to
$T$ the $n^2$-tuple formed by its entries (ordered in some way). Then there
exists a constant $C>0$ such that: $
H(T)\leq H\bigl(\phi (T)\bigr)\leq C H(T)$ for all $ T{\in}\Endk$.}
\endroster
Property (H4) follows at once from the definition of $\hop(T)$ and $H(T)$.
Of the remaining  properties, the only one which is not a straightforward consequence
of the corresponding properties of the local norms and the product formula is the inequality
$H\bigl(\phi (T)\bigr)\leq C H(T)$,
which is proven exactly as in lemma 2.1.
\remark{Remark} Note that H4 and lemma 2.1 implies that
$\hopf (T)<\infty$ for all $T{\in}\End{V}$, $V$ being a finite dimensional $K$-vector space and
$\F$ any adelic norm on it. 
\endremark
\vskip.2in
\head 3.  Comparison Results for Heights on $\Endk$ \endhead
The main result of this  section is theorem 3.3., which establishes a comparison result
between $H$ and 
$\hop$, which, beside of being interesting in its own right, 
will be used in the proof of theorem A. 
In fact our proof of  theorem A consists of two steps: first we prove the 
desired limit formula with  $H$ in place of $\hop$, then we show, by means of theorem 3.3, that this forces the limit formula to hold also 
for $\hop$.
\par 
Let  $\X\subseteq K^n$ be a subspace. For each $v{\in}\M{}{K}$ let $\X_v$ be the 
closure of $\X$ in $K_v^n$. Define  $\|\cdot\|_{\X_v}$, the {\it
seminorm relative to} $\X_v$, to
be $$ 
\nxv{\y}=\inf_{\x{\in}\X_v} \nv{\y-\x}.
$$
The global function associated to the local seminorms is
$$ 
\align d_{\X}:K^n\ba\X&\longrightarrow \R\\
\y&\mapsto\dx{\y}=\prod_{v{\in}\M{}{K}} \|\y\|_{\X_v}^{d_v}.\endalign
$$
The significance of $d_{\X}$ is explained by the following proposition
\proclaim{Proposition 3.1} {\it Suppose $T{\in}\Endk$ is not zero. Let
$\X=\ker{T}$. Then
$$
H(T)=\sup_{\y{\in}K^n\ba\X}
\frac{H\bigl(T(\y)\bigr)}{\dx{\y}}.\tag{3.1}
$$}
\endproclaim
\demo{Proof} 
$$
\nv{T}{=}
\sup_{\z{\in}K_v^n\ba\X_v}\{\sup_{\x{\in}\X_v}
\frac{\nv{T(\z)}}{\nv{\z-\x}}\}{=}\sup_{\z{\in}K_v^n\ba\X_v}
\{\frac{\nv{T(\z)}}{\underset{\x{\in}\X_v}\to\inf{\nv{\z-\x}}}\}{=}
\sup_{\z{\in}K_v^n\ba\X_v}
\frac{\nv{T(\z)}}{\nxv{\z}},
$$
hence
$$
\sup_{\y{\in}K^n\ba\X}\frac{
H\bigl(T(\y)\bigr)}{\dx{\y}}\leq \prod_{v{\in}\M{}{K}} 
\sup_{\z{\in}K_v^n\ba\X_v}\frac{\normdv{T(\z)}}{\nxdv{\z}}= H(T).
$$
To prove the reverse inequality,
let $\cs\supset\M{\infty}{K}$ be a finite set of places such that
$$
\nv{T}=1\qquad \forall v{\notin}\cs \qquad\text{and}\qquad
\frac{H\bigl(T(\z)\bigr)}{\dx{\z}}=\underset{v{\in}\cs}\to{\prod}\frac{\normdv{T(\z)}}
{\nxdv{\z}}\qquad\forall \z{\in}K^n\ba\X.\tag{3.2}
$$
The existence of $\cs$ is guaranteed by lemma 3.2 below.
Given $\epsilon >0$ let 
$\delta >0$ be such that $
\underset{v{\in}\cs}\to{\prod}\normdv{T} \leq\epsilon
+\underset{v{\in}\cs}\to{\prod}\bigl(\normdv{T}-\delta\bigr)
$.
By the weak approximation theorem we can find $\z{\in}K^n$ such that 
$$
\normdv{T}-\delta\leq\frac{\normdv{T(\z)}}{\normdv{\z}}
$$
for all $v{\in}\cs$. Taking the product over $\cs$ and using the equalities $(3.2)$ yields
$$
H(T)=\underset{v{\in}\cs}\to{\prod}\normdv{T} \leq\epsilon
+\underset{v{\in}\cs}\to{\prod}\bigl(\normdv{T}-\delta\bigr)\leq \epsilon+
\underset{v{\in}\cs}\to{\prod}\frac{\normdv{T(\z)}}{\normdv{\z}}
\leq\epsilon+ \frac{H\bigl(T(\z)\bigr)}{\dx{\z}}
$$
completing the proof of the lemma.$\fine$\enddemo

\proclaim{\bf Lemma 3.2} {\it Let $0\neq T{\in}\Mn\bigl(K\bigr)$ and set $\X=\ker{T}$.
Then there exists a finite set of places $\cs\supset\M{\infty}{K}$
\roster
\item"(a)" $\nv{T}=1\qquad \forall v{\notin}\cs$
\item"(b)" ${\displaystyle\frac{H\bigl(T(\z)\bigr)}{\dx{\z}}=\underset{v{\in}\cs}\to{\prod}\frac{\normdv{T(\z)}}
{\nxdv{\z}},\qquad \forall \z{\in}K^n\ba\X} $.
\endroster}
\endproclaim
\demo{Proof} If $\rank{T}=n$ the lemma reduces to the fact that $T$ belongs to 
$\operatorname{GL}_n\bigl(\Ov\bigr)$
for all but finitely many $v{\in}\M{}{K}$.
So we can assume $r<n$. To prove the lemma it suffices to exhibit a subspace $\Y\subset K^n$ of dimension $r=\rank{T}$ and a finite set of places
$\cs\supset\M{\infty}{K}$ such that for all $v{\notin}\cs$ we have :
\roster
\item"(a$^\prime$)" $\nxv{\y}= \nv{\y}\qquad\qquad\forall\y{\in}\Y_v$ .
\item"(b$^\prime$)" $
\nv{T(\y)}=\nxv{\y}\qquad\qquad\forall\y{\in}\Y_v$ .
\endroster
Suppose first that $\X$ is spanned by the last $n-r$ elements of the canonical basis of $K^n$. 
In this case we take $\Y$ to be the span of the first $r$ elements of the 
canonical basis of $K^n$, and properties (a$^\prime$) and (b$^\prime$) are trivially  verified. 
In general we proceed as follows: choose $S{\in}
\Mn\bigl(K\bigr)$, $S$ invertible, such that $\X$ is 
spanned by $S(\e_{r+1}),\ldots  , S(\e_{n})$.  Now note that
 (1) and (2) hold for $T\comp S$ and that
 $S$ is an isometry with respect to $\nv{{\cdot}}$ for all but finitely many $v{\in}\M{}{K}$. Let 
$\Y=S(<\e_1, \ldots , \e_r>)$ and $\cs$ be the finite set formed 
by those places of $K$
for which either $S$ is not a $\nv{{\cdot}}$-isometry, or one of conditions (a$^\prime$) and 
(b$^\prime$) does not hold
for $T\comp S$. It is straightforward to verify that $\Y$ and $\cs$
satisfies (a$^\prime$) and (b$\prime$).
$\fine$ \enddemo
\proclaim{\bf Theorem 3.3} {\it Let
$T{\in}\End{K^n}$, and set $\X=\ker{T}$. 
Then
\roster \item"(a)" If $T$ is invertible, then
$
\hop (T) =H(T).
$
\item"(b)" If $1<\rank{T}<n$, then there a constant $C(K,n)\geq 1$ depending only on $K$ and $n$.
$$
 \hop (T)\geq H(T){C(K,n)}^{-1} H(\X)^{-1} . 
$$
\item"(c)" If $\rank{T}=1$, then
$
\hop (T) = H(T)H(\X)^{-1}.
$
\endroster}
\endproclaim
The proof of theorem 3.3 is based upon the following lemma:
\proclaim{\bf Lemma 3.4} 
{\it \roster\runinitem"(a)" Let $\X\subset K^n$ be a subspace
such that  $1\leq\dim{K}{\X}< n-1$. Then 
$$
\inf_{\x{\in}\X}H(\y-\x)\leq C(K,n) \dx{\y} H(\X)
$$
 for all
$\y{\notin}\X$.
\item"(b)" If $\dim{K}{\X}= n-1$, then 
$
H(\X)=d_{\X}(\y)^{-1}.
$
\endroster}
\endproclaim 
\demo{Proof}
It follows from lemma 4 of \cite{12} that
$
H(\X)\dx{\y}=H\bigl(\langle \X ,\y\rangle\bigr).
$
If $\dim{K}{\X}=n-1$, then  $\langle \X,
\y \rangle=K^n$, proving (b) (recall that $H(K^n)=1$). Now 
assume $\dim{K}{\X}< n-1$.
By applying the $\ell^2$-version of Siegel's lemma (see
\cite{11}), we find a   a basis $\{\z_1, \ldots , \z_{\ell +1}\}$ of $\langle \X
,\y\rangle$ and a constant $C(n,K)$ depending on $n$ and $K$, but not on $\X$, such that
$$
\underset{i=1}\to{\overset{\ell +1}\to\prod} H(\z_i)\leq C(n, K) H\bigl(\langle
\X,\y\rangle\bigr)=C(n,K) H(\X)\dx{\y}
$$
Now the height of any $\y{\in}K^n$ is at least one and at least one of the $\z_i$'s 
has to be of the form $\z_i= \y-\x_i$, completing the proof of the proposition. $\fine$
\enddemo
\demo{Proof of theorem 3.3}  \newline \noindent (a) If $\rank{T}=n$ then $d_{\X}=H$ and so
(a) was proven in proposition 3.1
\newline \noindent
(b) Let $\X=\ker{T}$ and $C=C(n,K)$. By proposition 3.1 and lemma 3.4., we have
$$
H(T)=\sup_{\y{\notin}\X}
\frac{H\bigl(T(\y)\bigr)}{\dx{\y}} \leq C\cdot
H(\X)\cdot\left\{\sup_{\y{\in}K^n} \frac{H\bigl(T(\y)\bigr)}{\inf_{\x{\in}\X}H(\y-\x)}\right\}
=C\cdot
H(\X) \hop (T).
$$
\noindent (c) Let $\Y$ be the image of $T$. Then $\hop (T)=H(\Y)=H(T(\y))$ for
any $\y{\notin}\X$. By proposition 3.1 $H(T)=H(T(\y))d_{\X}(\y)^{-1}$, hence (c) follows from
 lemma 3.4. $\fine$
\enddemo
\vskip.2in
\head 4.  Proof of Theorems A and B\endhead
The main ingredients for the proof of theorem A, besides theorem 3.3, 
are the local Gelfand-Beurling formulae, which
we recall below:
\proclaim{Theorem 4.1}{\it Let $K$ be a number field and $v$ a place of $K$. Suppose  $S$ belongs to $\End{K_v^n}$. Then
$$\lim_{k\to\infty}\|S^k\|_v^{{}^{\ssize \frac1k}}=\rov{S}.\tag{4.1}
$$}
\endproclaim
\demo{\it Proof} If $v{\in}\M{\infty}{K}$, then (4.1)
is a special case of the Gelfand-Beurling formula for Banach algebras.
For a simultaneous proof of the real and complex case see \cite{3}. 
If $v{\in}\M{\infty}{K}$, then (4.1) is proven (in a more general setting)
in \cite{1, theorem 7.2.1.} A direct and elementary proof  is in
\cite{10, supplement 3, theorem 14}.
\enddemo
The next theorem completes the proof of theorem A of the introduction.
\proclaim{Theorem 4.2} {\it Let $T{\in}\Endk$. Then  $$
\sph  (T) = \lim_{k\to\infty} \bigl(\hop (T^k)\bigr)^{\frac1k}. $$
}
\endproclaim
\demo{\it Proof} We will first prove that
$$
\sph  (T) = \lim_{k\to\infty} \bigl(H(T^k)\bigr)^{\frac1k}.\tag{5.2}
$$
Then we will deduce the analogous formula for $\hop$
with the aid of
theorem 3.3.
Fix $T{\in}\Endk$ and let $\cs\subset\M{}{K}$ be  defined by
requiring that $v{\in}\cs$ if and only if either $\rov{T}\neq 1$ or
$\nv{T^k}\neq 1$ for some $k\geq 1$. Now the first condition is clearly verified only
for finitely many $v$'s. By lemma 3.2 the second condition is also verified only by
finitely many $v$'s. It follows that $\cs$ is a finite set.
Moreover:
$$
\sph (T)= \underset
{v{\in}\cs}\to{\prod}\rovdv{T_v};\qquad\qquad
H({T^k})=\underset {v{\in}\cs}\to{\prod} \normdv
{T^k_v}
\qquad\qquad\text{
for all } k\geq 1.
$$
Since $\cs$ is finite we can exchange the product with
the limit,  formula (4.2) follows from the local Gelfand-Beurling
formulae. In particular, this yields the theorem for $T$ invertible because in this
case $H(T)=\hop (T)$ by theorem 3.3.(a). Suppose now that $T$ is singular. 
Note that $\ker{T^k}=\ker{T^h}$ for all $h,k\geq n$. Let
$B=H(\ker{T^n})$. Then, by theorem 3.3.(b) and (c), we have
$$
\liminf_{k\to\infty}\hop (T^k)^{\frac1k}\geq \liminf_{k\to\infty}
\left(\frac{1}{C_KB}\right)^{\frac1k} H(T^k)^{\frac1k}
\geq\liminf_{k\to\infty}H(T^k)^{\frac1k}=\sph (T).
$$
On the other hand  $\hop (T)\leq
H(T)$, so $$
\limsup_{k\to\infty} \hop (T^k)^{\frac1k}\leq \limsup_{k\to\infty} H (T^k)^{\frac1k}
=\sph (T).\fine
$$
\enddemo
\par
Next we deal with  Northcott's finiteness theorem. As was previously 
mentioned, it
is indeed possible to have an infinite family of rank-one maps which are 
pairwise not homothetic and which all have height bounded by the same 
constant. Consider for example the family 
$\bigl\{[T_n], n{\in}{\Bbb Z}\bigr\}$,
where
$
T_n=\left(\smallmatrix 1&n\\1&n\endsmallmatrix\right)$. Then 
$\hop(T_n)=\sqrt{2}$ for all $n$, and $T_n$ is not homothetic to $T_m$, if $n\neq m$. Let
$$
\Omega_1^C\bigl(\P,B\bigr)=\left\{\cl{T}
 \st \rank{T}=1,\hop (T)\leq B\text{ and }
H\bigl(\ker{T}\bigr)\leq C\right\}.
$$
Then our finiteness results can be formally stated as follows:
\proclaim{\bf Theorem 4.3 } {\it Let $B\geq 1$. Then
\roster
\item"(a)"  The set $\Omega\bigl(\P,\hop ,B\bigr)
$
is finite.
\item"(b)" The set
$
\Omega_1^C\bigl(\P,B\bigr)
$ is finite for all $C\geq 1$.
\endroster}
\endproclaim
\demo{Proof} (a)  Northcott's finiteness theorem 
for projective space implies that we can choose
 $\zero=\x_0,\x_1, \ldots ,\x_m{\in}K^n$, in such a way 
that any $\y{\in}K^n$ with $H(\y)\leq B$ is a non-zero scalar multiple of one and only one of the
$\x_i$'s. Suppose $T=(\t_1\,\,\ldots \,\,\t_n)$ is such that $\hop (T)\leq B$. 
 Then
$\t_j=\lambda_j (T)\x_{i^T_j}$, with $\lambda_j\neq 0$ and $i^T_j{\in}\{0, \ldots , m\}$. The $n$-tuple 
$(i^T_1, \ldots ,i^T_n)$ depends only on the image of $T$ in $\P$. To prove (a)
it suffices to show that there are only finitely many elements of $\om$
having the same associated  $n$-tuple. Fix an $n$-tuple  $(i_1, \ldots ,i_n)$ which arises
as associated to some element $[T]=[(\t_1\,\,\ldots \,\,\t_n)]{\in}\om$.
Since $\rank{T}\geq 2$, there exists $h,k$ such that 
$\zero\neq\x_{i_h}\neq \x_{i_k}\neq\zero$, and
$$
H(\x_{i_h}+\lambda_h(T)^{-1}\lambda_{j}(T)\x_{i_j})\leq
 \sqrt{2} B;\quad H(\x_{i_k}+\lambda_k(T)^{-1}\lambda_{j}(T)\x_{i_j})
\leq \sqrt{2}B\tag{4.3}
$$
for all $j$'s for which $\x_{i_j}\neq \zero$.  Northcott's Theorem for projective spaces implies that given $\zero\neq\y$ and $\z$ linearly independent, there are only finitely many
values of $\lambda{\in}K^{\times}$ for which vectors of the form $\y+\lambda \z$ have bounded height. 
Combining this with the inequalities $(4.3)$, we find 
that the ratios $\lambda_h(T)^{-1}\lambda_{j}(T)$ can assume only finitely many values. 
Hence  
$(i_1, \ldots ,i_n)$ is associated only to finitely many $[T]{\in}\om$. 
\newline \noindent (b)
Theorem 3.3. implies that 
$$\Omega_1^C\bigl(\projkn, B)\subset\bigl\{T{\in}
\projkn\st H(T)\leq B\cdot C\bigr\}.
$$
But  the set on the right is  a finite set by (H5) and Northcott's 
finiteness theorem for projective spaces. $\fine$
\enddemo

\enddocument
\end